
\input amstex
\documentstyle{amsppt}
\magnification=1200
\pagewidth{6.5truein}
\pageheight{9.0truein}
\TagsOnRight
\NoRunningHeads
\footline={\ifnum\pageno>1 \hfil\folio\hfil\else\hfil\fi}
\topmatter
\title\nofrills
  Special Isogenies and Tensor Product Multiplicities
\endtitle
\author
  Shrawan Kumar and John R. Stembridge
\endauthor
\date 14 June 2007\enddate
\address
  S.K.: Dept. of Mathematics, University of North Carolina,
  Chapel Hill NC 27599--3250
\endaddress
\email shrawan\@email.unc.edu\endemail
\address
  J.S.: Dept. of Mathematics, University of Michigan,
  Ann Arbor MI 48109--1043
\endaddress
\email jrs\@umich.edu\endemail
\thanks
  The authors were supported by NSF and NSA grants respectively.
  We thank J.~Millson for raising questions which led to this work,
  and we thank G.~Prasad for helpful discussions.
\endthanks
\abstract
We show that any bijection between two root systems that preserves
angles (but not necessarily lengths) gives rise to inequalities
relating tensor product multiplicities for the corresponding complex
semisimple Lie groups (or Lie algebras). We explain the inequalities
in two ways: combinatorially, using Littelmann's Path Model,
and geometrically, using isogenies between algebraic groups defined
over an algebraically closed field of positive characteristic.
\endabstract
\endtopmatter
\document


\def\mysec#1{\bigbreak\centerline{\bf #1}\message{ * }\nopagebreak\par}
\def\subsec#1{\medbreak{\it\noindent #1.}\nopagebreak\par}
\def\myref#1{\medskip\item"{[{\bf #1}]}"} 
\def\pf{{\it Proof.\ }} 
\def\endpf{\hbox{\qed}\medbreak}

\catcode`\@=11
\def\cite#1{\relaxnext@
  \def\nextiii@##1,##2\end@{[{\bf##1},\,##2]}%
  \in@,{#1}\ifin@\def\next{\nextiii@#1\end@}\else
  \def\next{[{\bf#1}]}\fi\next}
\def\proclaimheadfont@{\smc}
\let\proclaim\undefined
\outer\def\proclaim{%
  \let\savedef@\proclaim \let\proclaim\relax
  \add@missing\endroster \add@missing\enddefinition
  \add@missing\endproclaim \envir@stack\endproclaim
 \def\proclaim##1{\restoredef@\proclaim
   \penaltyandskip@{-100}\medskipamount\varindent@
   \def\usualspace{{\proclaimheadfont@\enspace}}\indent\proclaimheadfont@
   \ignorespaces##1\unskip\frills@{.\enspace}%
  \sl\ignorespaces}%
 \nofrillscheck\proclaim}
\def\remarkheadfont@{\smc}
\def\remark{\let\savedef@\remark \let\remark\relax
  \add@missing\endroster \add@missing\endproclaim
  \envir@stack\endremark
  \def\remark##1{\restoredef@\remark
    \penaltyandskip@\z@\medskipamount
  {\def\usualspace{{\remarkheadfont@\enspace}}\indent%
  \varindent@\remarkheadfont@\ignorespaces##1\unskip%
  \frills@{.\enspace}}\rm
  \ignorespaces}\nofrillscheck\remark}
\def\endremark{\par\revert@envir\endremark
  \penaltyandskip@{55}\medskipamount}
\catcode`\@=\active


\define\om{\omega}
\define\al{\alpha}

\define\lam{\lambda}

\define\eps{\varepsilon}

\define\br#1{\langle #1\rangle}
\define\widebar#1{\overline{#1}}
\def\ul{^{\ell}}
\def\dl{_{\ell}}

\define\simto{\overset\sim\to\rightarrow}
\define\tosim{\hbox to 16pt{\hfil
  \rlap{\kern 3pt\raise -3pt\hbox{$\scriptstyle\sim$}}
  \kern -2pt $\rightarrow$\hfil}}

\define\ch{\operatorname{ch}}

\define\chr{\operatorname{char}}
\define\Span{\operatorname{Span}}
\define\SO{\operatorname{SO}}
\define\SP{\operatorname{Sp}}
\define\Spin{\operatorname{Spin}}

\define\bq{\Bbb Q}
\define\bbr{\Bbb R}
\define\bz{\Bbb Z}
\define\bc{\Bbb C}

\loadeusm
\define\cl{{\eusm L}}
\define\LS{{\eusm C}}
\define\tp{T}

\define\Lman{Littelmann}
\define\xtal{crystallographic}


{\openup 0.5\jot   

\mysec{1. Introduction}
Let $G=G(k)$ and $G'=G'(k)$ be connected, semisimple algebraic
groups over an algebraically closed field $k$ of
characteristic $p>0$, and let $f:G\to G'$ be an isogeny
(i.e., a surjective algebraic group homomorphism with a finite kernel).
For example (see~\cite{BT,\S3.3}), there is an isogeny
$\SO_{2\ell+1}(k)\to\SP_{2\ell}(k)$ in characteristic~2.

Fix a Borel subgroup $B$ of $G$ and $T\subset B$ a maximal torus.
Given that $f$ is an isogeny, it follows that $T':= f(T)$ is a
maximal torus in $G'$ and $B':=f(B)$ is a Borel subgroup of~$G'$.  
Letting $X(T)$ denote the group of characters of $T$
and $X(T)^+\subset X(T)$ the subset of dominant characters,
the map $f$ induces a homomorphism $f^*: X(T')\to X(T)$
that sends $X(T')^+$ into $X(T)^+$.

Now, let $G(\bc)$ denote the connected semisimple complex
algebraic group with the same root datum as $G(k)$.
In this paper, we show that for each isogeny $f:G(k)\to G'(k)$,
there is a family of inequalities relating tensor product
multiplicities for $G(\bc)$ and $G'(\bc)$. More explicitly,
for each $\lambda\in X(T)^+$, let $V(\lambda)$ denote the
irreducible representation of $G(\bc)$ with highest weight $\lambda$,
and for any sequence $\lambda_1,\dots,\lambda_n\in X(T)^+$, let
$$
[\lam_1,\dots,\lam_n]^{G(\bc)}:=
  \dim\bigl(V(\lam_1)\otimes\cdots\otimes V(\lam_n)\bigr)^{G(\bc)}\tag1.1
$$
denote the dimension of the $G(\bc)$-invariant subspace of the
corresponding tensor product. In these terms, we prove in
Theorem~3.3 below that for all $\lam'_1,\dots,\lam'_n \in X(T')^+$,
$$
[\lam'_1,\dots,\lam'_n]^{G'(\bc)}\leq
  [f^*(\lam'_1),\dots,f^*(\lam'_n)]^{G(\bc)}.\tag1.2
$$

For example, using the characteristic~2 isogeny
$f:\SO_{2\ell+1}(k)\to\SP_{2\ell}(k)$,
the standard coordinates one finds in the appendices of~\cite{B}
lead to identifications of both $X(T)$ and $X(T')$ with the
lattice $\bz^\ell$ so that the induced map $f^*:\bz^\ell\to\bz^\ell$
is the identity map (cf.~Section~4.A). Thus, as a special case of~(1.2),
we obtain the inequality
$$
[\lam_1,\dots,\lam_n]^{\SP_{2\ell}(\bc)}
  \leq [\lam_1,\dots,\lam_n]^{\SO_{2\ell +1}(\bc)}
$$
for all dominant $\SP_{2\ell}(\bc)$-weights $\lambda_1,\dots,\lambda_n$.
In Section~4, we will return to this inequality and three more of a
similar flavor, and deduce from two of the inequalities an interesting
saturation property for tensor products (see Remark~4.5(b)).

We give two proofs of~(1.2), one combinatorial and one geometric.

In the combinatorial approach,
we extract the essential necessary features of the induced
maps $f^*$ by introducing the notion of an {\it integer
renormalization} $\phi:X(T')_{\bbr}\simto X(T)_{\bbr}$,
where $X(T)_{\bbr}:=X(T)\otimes_\bz\bbr$. This is purely a root
system concept, not tied to the algebraic group context,
and mildly generalizes the {\it special isomorphisms} arising from
isogenies (see Definition~3.1 and Proposition~3.2 below).

We then prove a slight generalization of~(1.2) (see Theorem~2.1 and
Corollary~2.2), replacing $f^*$ with any integer renormalization $\phi$.
The proof is based on \Lman's Path Model for tensor product
multiplicity (see~\cite{L1} and~\cite{L2}). More specifically,
we use a variant of the Path Model (see~\cite{S}) in which the
objects are chains in the Bruhat ordering of various
Weyl group orbits, and the inequality is obtained by comparing
chains related by integer renormalizations.

In the geometric approach, we use an important result
due to Donkin asserting the existence of good filtrations
for products of global sections of homogeneous line bundles.
(It should be noted that
Donkin proved this result for almost all the cases barring a few
exceptions involving small primes~\cite{D}; the result was
subsequently proved uniformly by Mathieu for all primes~\cite{M}.)
This allows replacing~(1.2) with a cohomological statement that is
independent of the characteristic of the field
(including the characteristic 0 case), thereby enabling
us to deduce~(1.2) directly from the existence of an isogeny
in characteristic $p$.

\mysec{2. The combinatorial approach}
Let $R$ and $R'$ be finite \xtal{} root systems embedded
(respectively) in real Euclidean spaces $E$ and $E'$ with inner
products $\br{\cdot\,{,}\,\cdot}$ and $\br{\cdot\,{,}\,\cdot}'$.
If there is an isometry $\phi:E'\to E$ and a map $c:R\to\bz^{>0}$
such that
$$
R'=\{c(\alpha)\phi^{-1}(\alpha):\alpha\in R\},
$$
then we say that the pair $(\phi,c)$ defines an
{\it integer renormalization} from $R'$ to $R$.

We may embed the dual root system
$R^\vee=\{\alpha^\vee:\alpha\in R\}$ in $E$
(and similarly embed $(R')^\vee$ in $E'$) by setting
$\alpha^\vee:=2\alpha/\br{\alpha,\alpha}$. In these terms,
one easily checks that if $(\phi,c)$ is an integer renormalization
$R'\to R$ and $\alpha'=c(\alpha)\phi^{-1}(\alpha)\in R'$, then
$$
(\alpha')^\vee={1\over c(\alpha)}\phi^{-1}(\alpha^\vee).\tag2.1
$$
Thus, there is also an integer renormalization $(\phi^{-1},c')$ 
from $R^\vee$ to $(R')^\vee$, where $c'$ is obtained by
setting $c'((\alpha')^\vee)=c(\alpha)$. Furthermore, letting
$$
P(R):=\{\lambda\in\Span_{\Bbb R} R:\alpha\in R\Rightarrow
  \br{\lambda,\alpha^\vee}\in\bz\}
$$
denote the lattice of integral $R$-weights, (2.1) implies that
$$
\br{\lambda',(\alpha')^\vee}'\in\bz\ \Rightarrow\ 
  \br{\phi(\lambda'),\alpha^\vee}\in c(\alpha)\bz,
$$
and thus $\phi$ injects $P(R')$ into $P(R)$.

For example, if we set $E'=E$ and take $\phi$ to be the identity map,
then a trivial renormalization may be obtained by
setting $c(\alpha)=c$ for some fixed positive integer~$c$.
Here, $R'=cR$ is isomorphic to $R$. 

Less trivially, assume $R$ has both long and short roots,
with the squared ratio of long root lengths to short root
lengths being $r$. (This is possible only if $r=2$ or $3$.)
Again keeping $E'=E$ and $\phi$ the identity map, setting
$$
c(\alpha)=\cases r&\text{if $\alpha$ is short},\\
  1&\text{if $\alpha$ is long}\endcases\tag2.2
$$
also yields an integer renormalization.
In this case, $R'$ is isomorphic to $R^\vee$.

Given a root system $R$, fix a choice of positive roots $R^+\subset R$,
and let $P(R)^+\subset P(R)$ denote the resulting semigroup of
dominant weights. If $(\phi,c)$ is an integer renormalization from
$R'$ to~$R$, we require that a compatible choice of positive
roots is made in $R'$ and $R$ so that
$(R')^+=\{c(\alpha)\phi^{-1}(\alpha):\alpha\in R^+\}$.
With this choice, $\phi$ injects  $P(R')^+$ into  $P(R)^+$.

Now let $G(\bc)$ be a connected semisimple algebraic group over
$\bc$ with root system $R$. For convenience (see Remark~2.3 below),
we assume $G(\bc)$ is simply connected, so that
the character lattice $X(T)$ may be naturally identified with $P(R)$.
Thus for each $\lambda\in P(R)^+$, there is an irreducible
$G(\bc)$-representation $V(\lambda)$ of highest weight $\lambda$,
and we may define
$$
m_{G(\bc)}(\lambda;\mu,\nu):=\text{multiplicity of $V(\lambda)$
  in $V(\mu)\otimes V(\nu)$}
$$
for all $\lambda,\mu,\nu\in P(R)^+$.

\proclaim{Theorem 2.1}
If $G(\bc)$ and $G'(\bc)$ are two connected, simply connected,
semisimple algebraic groups over $\bc$ with root systems $R$ and $R'$,
and $(\phi,c)$ is an integer renormalization from $R'$ to $R$,
then for all $\lambda',\mu',\nu'\in P(R')^+$, we have
$$
m_{G'(\bc)}(\lambda';\mu',\nu')
  \le m_{G(\bc)}(\phi(\lambda');\phi(\mu'),\phi(\nu')).
$$
\endproclaim

\pf \Lman's Path Model (see~\cite{L1} and~\cite{L2}) provides an
explicit combinatorial description for $m_{G(\bc)}(\lambda;\mu,\nu)$.
Following the approach in Section~8 of~\cite{S}, the objects in this
model may be viewed as chains in the Bruhat orderings of various Weyl
group orbits.

More explicitly, let $W(R)=\br{\sigma_\alpha:\alpha\in R^+}\subset GL(E)$
denote the Weyl group generated by the reflections
$\sigma_\alpha:\lambda\mapsto\lambda-\br{\lambda,\alpha^\vee}\alpha$.
Given $\mu\in P(R)^+$, let `$<$' denote the Bruhat ordering
of the $W(R)$-orbit of~$\mu$; i.e., the transitive closure of the
relations $\sigma_\alpha(\nu)<\nu$ for all $\nu\in W(R)\mu$ and all
$\alpha\in R^+$ such that $\nu-\sigma_\alpha(\nu)$ is a {\it positive}
multiple of $\alpha$. We write $\nu_1\lessdot\nu_2$ when
$\nu_2$ covers $\nu_1$ in this order; i.e., $\nu_1<\nu_2$ and there
does not exist $\nu\in W(R)\mu$ such that $\nu_1<\nu<\nu_2$.

For each $b\in\bq^{>0}$, the $b$-Bruhat ordering of $W(R)\mu$ is
the partial order `$<_b$' generated by imposing the relations
$\sigma_\alpha(\nu)<_b\nu$ ($\nu\in W(R)\mu$, $\alpha\in R^+$)
whenever $\sigma_\alpha(\nu)\lessdot\nu$ and
$b(\nu-\sigma_\alpha(\nu))$ is an {\it integer} multiple of $\alpha$.
If $b\in\bz^{>0}$, we recover the usual Bruhat ordering.

A {\it Lakshmibai-Seshadri $R$-chain $C$ of type} $\mu$ is
defined to be a sequence
$$
\mu_0<_{b_1}\mu_1<_{b_2}\cdots<_{b_\ell}\mu_\ell,\tag2.3
$$
where $0=b_0<b_1<\dots<b_{\ell+1}=1$ ($b_i\in \bq$),
$\mu_0,\dots,\mu_{\ell}\in W(R)\mu$, and $\ell\ge0$.

We let $\LS(R,\mu)$ denote the set of all such chains.

Given $C$ as in~(2.3), define a sequence
$0=\delta_0(C),\delta_1(C),\dots,\delta_{\ell+1}(C)\in E$ by setting
$$
\delta_t(C)=\sum_{j=1}^t(b_j-b_{j-1})\mu_{j-1}
  =b_t\mu_{t-1}-\sum_{j=1}^{t-1}b_j(\mu_j-\mu_{j-1}).
$$
It is clear from the definition that $\nu_1<_b\nu_2$ implies
that $b(\nu_2-\nu_1)$ is an integral weight.
Hence, from the second expression for $\delta_t(C)$
(and the fact that $b_{\ell+1}=1$) we see that
$\omega(C):=\delta_{\ell+1}(C)$ is an integral weight.

The {\it depth} of $C$ is the unique integral weight $\delta(C)$
such that for each simple root $\alpha$,
$$
\br{\delta(C),\alpha^\vee}
  =\min_{0\le t\le \ell+1}\br{\delta_t(C),\alpha^\vee}.\tag2.4
$$
However, we should note that it is not immediately clear from this
definition that $\delta(C)$ is integral, rather than merely rational.
(For a proof, see Section~8 of \cite{S}.)

The description of tensor product multiplicities provided by
\Lman's Path Model may be expressed in the above terms as follows:
$$
m_{G(\bc)}(\lambda;\mu,\nu)=\left|\{C\in \LS(R,\mu):
  \nu+\delta(C)\in P(R)^+,\ \nu+\omega(C)=\lambda\}\right|.
$$
For a proof, see Theorem~8.3 of~\cite{S}
(cf.~also the Decomposition Rule in~\cite{L1}).

Now fix a dominant weight $\mu'\in P(R')^+$ and set $\mu=\phi(\mu')$.
Note that $\phi$ induces a Weyl group isomorphism $W(R')\to W(R)$,
and hence also an isomorphism between the Bruhat ordering $<$ 
of the $W(R')$-orbit of $\mu'$ and the $W(R)$-orbit of $\mu$.
Moreover, for all $\alpha'=c(\alpha)\phi^{-1}(\alpha)\in (R')^+$ 
(for $\alpha\in R^+$) and $\nu'\in W(R')\mu'$,
if $\sigma_{\alpha'}(\nu')<_b\nu'$ is a defining relation
of the $b$-Bruhat ordering of $W(R')\mu'$, then
$b(\nu'-\sigma_{\alpha'}(\nu'))$ is an integer multiple of $\alpha'$,
and hence $b(\nu-\sigma_\alpha(\nu))$ is an integer multiple
of $c(\alpha)\alpha$, where $\nu=\phi(\nu')$.
Thus, $\sigma_\alpha(\nu)<_b\nu$ is also a relation of the $b$-Bruhat
ordering of $W(R)\mu$. It follows that $\phi$ induces an
injection of $\LS(R',\mu')$ into $\LS(R,\mu)$.

Note that this injection is compatible with the weight
sequence $\delta_t(\cdot)$ in the sense that
$\delta_t(\phi(C'))=\phi(\delta_t(C'))$ for
all $C'\in\LS(R',\mu')$. Moreover, $\phi$ transforms the simple
roots of $R'$ into positive multiples of the simple roots of $R$,
so~(2.1) and~(2.4) also imply $\delta(\phi(C'))=\phi(\delta(C'))$.
Finally, since $\phi$ identifies the dominant chambers in $E'$
and $E$, it follows that for each $\lam',\mu',\nu'\in P(R')^+$,
$\phi$ injects the chains $C'\in\LS(R',\mu')$
that are counted by $m_{G'(\bc)}(\lam';\mu',\nu')$
(i.e., $\nu'+\delta(C')\in P(R')^+$, $\nu'+\omega(C')=\lam'$)
into those counted by
$m_{G(\bc)}(\phi(\lam');\phi(\mu'),\phi(\nu'))$.\endpf

In terms of invariant multiplicity (recall~(1.1)), one knows that
$$
[\lam,\mu,\nu]^{G(\bc)}=m_{G(\bc)}(\lam^*;\mu,\nu),
$$
where $\lam^*$ denotes the highest weight of $V(\lam)^*$,
and more generally,
$$
[\lam_1,\dots,\lam_n]^{G(\bc)}=\sum_{\mu\in P(R)^+}
  m_{G(\bc)}(\mu;\lam_1,\lam_2)\cdot[\mu,\lam_3,\dots,\lam_n]^{G(\bc)}.
$$
Comparing this with the corresponding expansion for $G'(\bc)$,
an induction on $n$ yields

\proclaim{Corollary 2.2}
If $G(\bc)$, $G'(\bc)$, $R$, $R'$ and $(\phi,c)$ are as given in
Theorem~2.1, then for all $\lam'_1,\dots,\lam'_n\in P(R')^+$, we have
$$
[\lam'_1,\dots,\lam'_n]^{G'(\bc)}
  \le[\phi(\lam'_1),\dots,\phi(\lam'_n)]^{G(\bc)}.
$$
\endproclaim

\remark{Remark 2.3}
If we do not assume that $G(\bc)$ and $G'(\bc)$ are simply connected,
then the character lattices $X(T)$ and $X(T')$ are naturally
identified with sublattices of $P(R)$ and $P(R')$,
and all tensor product multiplicities for these groups can be
expressed in terms of tensor products for their simply connected
covering groups. Thus Theorem~2.1 and Corollary~2.2 remain valid in
this more general setting provided that we add the hypothesis that
$\phi$ maps $X(T')$ into $X(T)$. In particular, no extra hypothesis
is needed if we assume only that $G(\bc)$ is simply connected.
\endremark

\mysec{3. The geometric approach}
Returning to the setting of the introduction, let $G=G(k)$ and
$G'=G'(k)$ be connected, semisimple algebraic groups over an
algebraically closed field $k$ of characteristic $p>0$,
and let $f:G\to G'$ be an isogeny.
Fix a Borel subgroup $B$ of $G$ and $T\subset B$ a maximal torus,
and let $B'=f(B)$ and $T'=f(T)$ be the corresponding groups in $G'$.

Recalling that the map $f$ induces a
homomorphism $f^*: X(T')\to X(T)$,
we may extend $f^*$ to an isomorphism
$f^*_{\bbr}:X(T')_{\bbr}\simto X(T)_{\bbr}$,
where (as in the introduction) $X(T)_{\bbr}:=X(T)\otimes_{\bz}\bbr$.

Letting $R=R(G,T)$ denote the root system of $G$ with
respect to $T$ and similarly $R'= R(G',T')$, we recall the
following from Expos\'{e} n$^\circ$ 18, Definition~1 of~\cite{C}.

\remark{Definition 3.1}
An isomorphism $\phi : X(T')_{\bbr} \to X(T)_{\bbr}$ is called
{\it special} if
\roster
\item"{(a)}" $\phi(X(T'))\subset X(T)$, and
\item"{(b)}" there exist integers $d(\al)\ge0$ such that
$$
R'=\bigl\{p^{d(\al)}\phi^{-1}(\al):\al\in R\bigr\}.
$$
\endroster
\endremark

\vskip-\medskipamount
For any isogeny $f$ as above, the induced map $f^*_\bbr$ is a
special isomorphism. Conversely, for any special isomorphism
$\phi:X(T')_{\bbr}\to X(T)_{\bbr}$,
there exists an isogeny $f: G\to G'$ with $f^*_{\bbr}=\phi$
(cf.~Expos\'{e} n$^\circ$ 23, \S{3}, Th\'{e}or\`{e}me~1 of \cite{C}).

 To relate this to the setting of Section~2, one should identify
$X(T)_\bbr$ and $X(T')_\bbr$ with the spaces $E$ and $E'$,
and (cf.~Remark~2.3) $X(T)$ and $X(T')$ correspond to
sublattices of $P(R)$ and $P(R')$.
Moreover, as a consequence of Expos\'{e} n$^\circ$ 18,
Proposition~4 of~\cite{C}, a special isomorphism $\phi$ induces
an isomorphism between the Weyl groups $W(R')$ and $W(R)$,
and under this identification of the Weyl groups,
the map $\phi$ is Weyl group equivariant.
This allows one to impose Weyl group-compatible Euclidean metrics
on $E$ and $E'$ so that $\phi$ is an isometry, and thus yields
the following.

\proclaim{Proposition 3.2}
A special isomorphism $\phi:X(T')_{\bbr}\to X(T)_{\bbr}$ is a integer
renormalization from $R'$ to $R$.
\endproclaim

In the following result, we deduce the tensor product
inequalities in~(1.2) directly from the existence of an isogeny.

\proclaim{Theorem 3.3}
If $f: G\to G'$ is an isogeny of connected semisimple algebraic
groups over an algebraically closed field $k$ of characteristic $p>0$,
then for all $\lam'_1, \dots, \lam'_n \in X(T')^+$,
$$
  [\lam'_1,\dots,\lam'_n]^{G'(\bc)}\leq
  [f^*(\lam'_1),\dots,f^*(\lam'_n)]^{G(\bc)},
$$
where $G(\bc)$ is the connected semisimple complex algebraic group 
with the same root datum as that of $G(k)$ and similarly for $G'(\bc)$.
\endproclaim

\pf The map $f$ clearly induces a surjective morphism (of varieties)
$\bar{f}: X_n\to X'_n$, where $X_n := (G/B)^{\times n}$.
Consider the dominant line bundle
$\cl(\lam'_1)\boxtimes \cdots \boxtimes\cl (\lam'_n)$ on $X'_n$,
where $\cl (\lam'_i)$ is the homogeneous line bundle
$G'\times_{B'} k_{-\lam'_i}$ on $G'/B'$ corresponding to the
character $(\lam_i^\prime)^{ -1}$ of $B'$. In these terms,
the pull-back line bundle on $X_n$ is the homogeneous line bundle
$$
\cl (\lam_1)\boxtimes\cdots\boxtimes \cl (\lam_n),
$$
where $\lam_i:=f^*(\lam'_i)$. Thus, we get an injective map
$$
\bar{f}^*: H^0\bigl( X'_n, \cl (\lam'_1)\boxtimes \cdots \boxtimes \cl
  (\lam'_n)\bigr) \hookrightarrow H^0\bigl( X_n, \cl (\lam_1)\boxtimes
  \cdots \boxtimes \cl (\lam_n)\bigr) .
$$
Since the map $\bar{f}$ is $f$-equivariant under the diagonal action
of $G$ on $X_n$ and $G'$ on $X'_n$, the injection $\bar{f}^*$
induces an injection (still denoted by)
$$
\bar{f}^*: H^0\bigl( X'_n, \cl (\lam'_1)\boxtimes\cdots\boxtimes
  \cl (\lam'_n)\bigr)^{G'} \hookrightarrow
H^0\bigl(X_n,\cl(\lam_1)\boxtimes
  \cdots\boxtimes\cl(\lam_n)\bigr)^G.\tag3.1
$$
However, we have of course
$$
H^0\bigl( X_n, \cl (\lam_1) \boxtimes \cdots \boxtimes \cl(\lam_n)\bigr)
  \cong H^0\bigl( G/B, \cl (\lam_1)\bigr) \otimes \cdots \otimes
  H^0\bigl( G/B, \cl (\lam_n)\bigr).
$$
By Corollary~4.2.14 of~\cite{BrK}, the above module
$$
M := H^0\bigl( G/B, \cl (\lam_1)\bigr) \otimes\cdots\otimes
  H^0\bigl( G/B, \cl (\lam_n)\bigr)
$$
admits a good filtration.  Hence, by Theorem~4.2.7,
identity (4.2.1.3) and Proposition 4.2.3(c) of~\cite{BrK},
its $T$-character is
$$
\ch M = \sum_{\lam\in X(T)^+}\dim
  \bigl(H^0\bigl(G/B,\cl(\lam)\bigr)\otimes M\bigr)^G\cdot\ch(V_k(\lam)),
$$
where, for any $\lambda\in X(T)^+$, $V_k(\lambda):=H^0(G/B, \cl (\lam))^*$ 
is the Weyl module with highest weight $\lam$. Recall that,
by the Borel-Weil Theorem, 
$$
H^0\bigl( G(\bc)/B(\bc), \cl_\bc(\lam)\bigr)\simeq V(\lam)^*,
$$
where (as earlier) $ V(\lam)$ is the (complex) irreducible $G(\bc)$-module
with highest weight $\lambda$ and $ \cl_\bc(\lam)$ is the homogeneous
line bundle on $ G(\bc)/B(\bc)$ corresponding to the
character $\lambda^{-1}$ of $B(\bc)$. Moreover, as is well-known, 
$$
\ch V_k(\lam)=\ch V(\lam).
$$
(This follows from the vanishing of the cohomology
$H^i(G/B, \cl (\lambda))$ for all $i>0$.)

But clearly, $\ch M = \ch (V_k(\lam_1)^*)\cdots\ch (V_k(\lam_n)^*)$;
in particular, it is independent of the characteristic of the field
(including characteristic~0). Moreover, since
$\{\ch V(\lam )\}_{\lam\in X(T)^+}$ are $\bz$-linearly independent
as elements of the group ring of $X(T)$, we deduce that
$$
\dim\bigl( H^0\bigl( G/B, \cl (\lam )\bigr)\otimes M\bigr)^G
$$
is independent of the characteristic of the base field
for all $\lam\in X(T)^+$.
Taking $\lam=0$, we obtain that $\dim M^G$ is independent of the
characteristic. Observe next that (3.1) implies
$$
\dim (M' )^{G'} \leq \dim M^G,\tag3.2
$$
where
$$
M':= H^0\bigl( G'/B', \cl (\lam'_1)\bigr) \otimes\cdots\otimes
  H^0\bigl( G'/B', \cl (\lam'_n)\bigr).
$$
Thus~(3.2) implies
$$\align
  \dim \bigl(&V(\lam'_1) \otimes\cdots\otimes V(\lam'_n)\bigr)^{G'(\bc)}
    =\dim \bigl( V(\lam'_1)^* \otimes\cdots\otimes
  V(\lam'_n)^*\bigr)^{G'(\bc)}\\
    &\leq \dim \bigl( V(\lam_1)^* \otimes\cdots\otimes
  V(\lam_n)^*\bigr)^{G(\bc)} = \dim ( V(\lam_1)\otimes\cdots\otimes
  V(\lam_n))^{G(\bc)},
\endalign$$
and the theorem follows.\endpf

\mysec{4. Special cases}
An isogeny $f:G\to G'$ as in Section 3 for a simple $G$ is called
{\it special} if $d(\al)=0$ for some $\al\in R(G,T)$, where $d(\al)$
is as in Definition~3.1; it is {\it central} if $d(\al)=0$ for
all $\al\in R(G,T)$. A complete list of special non-central
isogenies may be found in~\cite{BT, \S3.3}. In the following,
we list the resulting tensor product inequalities implied by~(1.2).

\subsec{A. From $\SO$ to $\SP$}
Let $\chr{k}=2$, $\ell\ge 2$, $G$ the adjoint group of type $B\dl$
(i.e., $G=\SO_{2\ell+1}(k)$), and $G'$ the simply connected group
of type $C\dl$ (i.e., $G'=\SP_{2\ell}(k)$).
Following the notation from the appendices of~\cite{B},
we identify $X(T)=\bigoplus\ul_{i=1}\bz\eps_i$
and $X(T')=\bigoplus\ul_{i=1}\bz\eps_i$.
In these terms, the identity map
is a special isomorphism $X(T')_\bbr\to X(T)_\bbr$
giving rise to an isogeny
$$
f:\SO_{2\ell+1}(k)\to \SP_{2\ell}(k).
$$
Moreover, it identifies the dominant weights in $X(T')$ and $X(T)$
with respect to the choice of positive roots as in {\it loc.~cit.}
Since $f^*$ acts as the identity map, both of  Corollary~2.2 and
Theorem~3.3 specialize as follows.

\proclaim{Theorem 4.1}
If $\lam_1,\dots,\lam_n$ are dominant weights
for $\SP_{2\ell}(\bc)$, then
$$
[\lam_1,\dots ,\lam_n]^{\SP_{2\ell}(\bc)}
  \leq [\lam_1,\dots,\lam_n]^{\SO_{2\ell +1}(\bc)}.
$$
\endproclaim

\subsec{B. From $\SP$ to $\Spin$}
Let $\chr{k}=2$, $\ell\ge 2$, $G=\SP_{2\ell}(k)$, and $G'$ the
simply connected group of type $B\dl$ (i.e., $G'=\Spin_{2\ell+1}(k)$).
In this case, the coordinates from~\cite{B} identify
$X(T)=\bigoplus\ul_{i=1}\bz\eps_i$ as above,
and $X(T')=\{\sum_{i=1}^\ell a_i\eps_i:
  a_i\pm a_j\in\bz\text{ for all $i,j$}\}$.
In these terms, the map $\mu\mapsto2\mu$ defines a special
isomorphism $X(T')_\bbr\to X(T)_\bbr$ inducing an isogeny
$$
f:\SP_{2\ell}(k)\to \Spin_{2\ell+1}(k).
$$
Here, Corollary~2.2 and Theorem~3.3 specialize as follows.

\proclaim{Theorem 4.2}
If $\lam_1,\dots,\lam_n$ are dominant weights
for $\Spin_{2\ell+1}(\bc)$, then
$$
[\lam_1,\dots,\lam_n]^{\Spin_{2\ell +1}(\bc)}
  \leq [2\lam_1,\dots,2\lam_n]^{\SP_{2\ell}(\bc)}.
$$
\endproclaim

\subsec{C. The case of $F_4$}
Let $G=G'$ be of type $F_4$ and $\chr{k}=2$.
In this case, the simple roots generate $X(T)$.
Numbering them $\alpha_1,\alpha_2,\alpha_3,\alpha_4$ as in~\cite{B},
we have that $\alpha_1$ and $\alpha_2$ are long and
$\alpha_3$ and $\alpha_4$ are short. Moreover, there is a
special isomorphism $\phi:X(T)_\bbr\to X(T)_\bbr$ such that
$$
\phi(\al_1)=2\al_4,\ \phi(\al_2)=2\al_3,
  \ \phi(\al_3)=\al_2,\ \phi(\al_4)=\al_1.
$$
In terms of the corresponding fundamental weights
$\om_1,\om_2,\om_3,\om_4$, one has
$$
\phi(\om_1)=2\om_4,\ \phi(\om_2)=2\om_3,
  \ \phi(\om_3)=\om_2,\ \phi(\om_4)=\om_1.
$$
Thus, we obtain

\proclaim{Theorem 4.3}
If $\lam_1,\dots,\lam_n$ are dominant weights for $F_4(\bc)$, then
$$
[\lam_1,\dots,\lam_n]^{F_4(\bc)}
  \leq [\phi(\lam_1),\dots,\phi(\lam_n)]^{F_4(\bc)},
$$
where $\phi(a\om_1+b\om_2+c\om_3+d\om_4):=d\om_1+c\om_2+2b\om_3+2a\om_4$.
\endproclaim

\subsec{D. The case of $G_2$}
Let $G=G'$ be of type $G_2$ and $\chr{k}=3$.
Letting $\alpha_1$ and $\alpha_2$ denote the simple roots,
with $\alpha_1$ short and $\alpha_2$ long, there is a
special isomorphism $\phi:X(T)_\bbr\to X(T)_\bbr$ such that
$$
\phi(\al_1) = \al_2,\ \phi(\al_2) = 3\al_1.
$$
In terms of the corresponding fundamental weights 
$\om_1$ and $\om_2$, one has
$$
\phi(\om_1)=\om_2,\ \phi(\om_2)=3\om_1.
$$
Thus, we obtain

\proclaim{Theorem 4.4}
If $\lam_1,\dots,\lam_n$ are dominant weights for $G_2(\bc)$, then
$$
[\lam_1,\dots ,\lam_n]^{G_2(\bc)}
  \leq [\phi(\lam_1),\dots,\phi(\lam_n)]^{G_2(\bc)},
$$
where $\phi(a\om_1+b\om_2):=3b\om_1+a\om_2$.
\endproclaim

\remark{Remark 4.5}
(a) Any nonspecial isogeny is the composition of a power of the 
Frobenius homomorphism with a special isogeny. Now, Theorem~3.3
for the Frobenius homomorphism yields only the 
well-known (and easy to prove) inequality
$$
[\lam_1,\dots ,\lam_n]^{G(\bc)}
  \leq [p\lam_1,\dots, p\lam_n]^{G(\bc)}.
$$
Similarly, central isogenies do not yield any new inequalities.

(b) Because of its connection with the eigenvalue problem,
there has been considerable renewed interest in characterizing
nonvanishing tensor product multiplicities, and more specifically
(e.g., see~\cite{KM} and~\cite{BeK} and the references therein),
in the semigroups
$$
\tp_n(G):=\{(\lambda_1,\dots,\lambda_n)\in(X(T)^+)^n:
  [\lambda_1,\dots,\lambda_n]^{G(\bc)}>0\},
$$
and their $\bq$-saturated analogues
$$
\widebar{\tp}_n(G):=\{(\lambda_1,\dots,\lambda_n)\in(X(T)_\bq^+)^n:
  (N\lambda_1,\dots,N\lambda_n)\in\tp_n(G)\text{ for some $N>0$}\},
$$
where $X(T)_\bq^+$ denotes the dominant part of $X(T)\otimes_\bz\bq$.
Although these semigroups are difficult to describe explicitly
in general, it is interesting to note that an immediate corollary
of Theorems~4.1 and~4.2 is that
$\widebar{\tp}_n(\SP_{2\ell})=\widebar{\tp}_n(\Spin_{2\ell+1})$.
\endremark

\par}

\mysec{References}
{\frenchspacing\raggedbottom
\roster
\myref{BT}
  A. Borel and J. Tits, Homomorphismes ``abstraits'' de groupes
  alg\'ebriques simples, {\it Ann. of Math.} {\bf 97} (1973), 499--571.
\myref{B}
  N. Bourbaki, {\it Groupes et Alg\'{e}bres de Lie}, Ch. 4--6,
  Masson, Paris, 1981.
\myref{BeK}
  P. Belkale and S. Kumar, Eigenvalue problem and a new product in 
  cohomology of flag varieties, {\it Invent. Math.} 
  {\bf 166} (2006), 185--228. 
\myref{BrK}
  M. Brion and S. Kumar, {\it Frobenius Splitting Methods
  in Geometry and Representation Theory}, Birkh\"auser, 2005.
\myref{C}
  C. Chevalley, {\it S\'{e}minaire C. Chevalley 1956--58},
  Classification des Groupes de Lie Alg\'ebriques, vol. 1, Ecole
  Normale Sup\'erieure, 1958.
\myref{D}
  S. Donkin, {\it Rational Representations of Algebraic Groups},
  Lecture Notes in Math. {\bf 1140}, Springer-Verlag, 1985.
\myref{KM}
  M. Kapovich and J. J. Millson, Structure of the tensor product 
  semigroup, {\it Asian J. Math.} {\bf 10} (2006), 493--539.
\myref{L1}
  P. Littelmann, A Littlewood-Richardson rule for symmetrizable
  Kac-Moody algebras, {\it Invent. Math.} {\bf 116} (1994), 329--346.
\myref{L2}
  P. Littelmann, Paths and root operators in representation theory,
  {\it Ann. of Math.} {\bf 142} (1995), 499--525.
\myref{M}
  O. Mathieu, Filtrations of $G$-modules,
  {\it Ann. Sci. \'Ec. Norm. Sup\'er.} {\bf 23} (1990), 625--644.       
\myref{S}
  J. R. Stembridge, Combinatorial models for Weyl characters,
  {\it Advances in Math.} {\bf 168} (2002), 96--131.
\endroster\par}
\enddocument